\documentclass{amsart}
\input{amssymb.sty}
\title{Weak Amenability of Hyperbolic Groups}
\author{Narutaka OZAWA}
\address{Department of Mathematical Sciences,
University of Tokyo, Komaba, 153-8914\\ \indent 
Department of Mathematics, UCLA, Los Angeles, CA 90095-1555}
\thanks{Supported by Sloan Foundation and NSF grant}
\email{narutaka@ms.u-tokyo.ac.jp}

\newtheorem{thm}{Theorem}
\newtheorem{lem}[thm]{Lemma}
\newtheorem{prop}[thm]{Proposition}
\theoremstyle{definition}
\newtheorem*{defn}{Definition}

\newcommand{\C}{{\mathbb C}}
\newcommand{\N}{{\mathbb N}}
\newcommand{\Z}{{\mathbb Z}}
\newcommand{\B}{{\mathbb B}}
\newcommand{\D}{{\mathbb D}}
\newcommand{\NN}{{\mathfrak N}}
\newcommand{\PP}{{\mathfrak P}}
\newcommand{\G}{\Gamma}

\newcommand{\p}{\varphi}
\newcommand{\pathp}{{\mathfrak p}}
\newcommand{\hh}{{\mathcal H}}

\newcommand{\cb}{\mathrm{cb}}
\newcommand{\ip}[1]{\langle#1\rangle} 
\begin{document}
\begin{abstract}
We prove that hyperbolic groups are weakly amenable. 
This partially extends the result of Cowling and Haagerup 
showing that lattices in simple Lie groups of real rank one
are weakly amenable. 
We take a combinatorial approach in the spirit of 
Haagerup and prove that for the word length 
distance $d$ of a hyperbolic group, 
the Schur multipliers associated with the kernel 
$r^d$ have uniformly bounded norms for $0<r<1$. 
We then combine this with a Bo\.zejko-Picardello type inequality 
to obtain weak amenability. 
\end{abstract}
\maketitle
\section{Introduction}
The notion of weak amenability for groups
was introduced by Cowling and Haagerup \cite{ch}.
(It has almost nothing to do with the notion of weak amenability
for Banach algebras.)
We use the following equivalent form of the definition.
See Section~\ref{sec:schur} and \cite{bo,ch,pisier} for more information.
\begin{defn}
A countable discrete group $\G$ is said to be 
\emph{weakly amenable with constant $C$} if there exists 
a sequence of finitely supported functions $\p_n$ on $\G$ 
such that $\p_n\to 1$ pointwise and $\sup_n\|\p_n\|_{\cb}\le C$, 
where $\|\p\|_{\cb}$ denotes the (completely bounded) norm 
of the Schur multiplier on $\B(\ell_2\G)$ associated with 
$(x,y)\mapsto \p(x^{-1}y)$.
\end{defn}

In the pioneering paper \cite{haagerup}, Haagerup 
proved that the group C$^*$-algebra of a free group has 
a very interesting approximation property. 
Among other things, he proved that the graph distance 
$d$ on a tree $\G$ is conditionally negatively definite; 
in particular, the Schur multiplier on $\B(\ell_2\G)$ 
associated with the kernel $r^d$ has (completely bounded) 
norm one for every $0<r<1$. 
For information of Schur multipliers and completely 
bounded maps, see Section~\ref{sec:schur} and \cite{bo,ch,pisier}. 
Bo\.zejko and Picardello \cite{bp} proved that 
the Schur multiplier associated 
with the characteristic function of the subset 
$\{ (x,y) : d(x,y)=n\}$ 
has (completely bounded) norm at most $2(n+1)$. 
These two results together imply that a group acting properly on 
a tree is weakly amenable with constant one. 
Recently, this result was extended to the case of 
finite-dimensional CAT(0) cube complexes 
by Guentner and Higson \cite{gh}. See also \cite{mizuta}. 
Cowling and Haagerup \cite{dch,cowling,ch} proved that 
lattices in simple Lie groups of real rank one
are weakly amenable and computed explicitly the associated constants. 
It is then natural to explore this property for hyperbolic groups 
in the sense of Gromov \cite{gdh,gromov}. 
We prove that hyperbolic groups are weakly amenable, without giving 
estimates of the associated constants. 
The results and proofs are inspired by and partially generalize 
those of Haagerup \cite{haagerup}, Pytlik-Szwarc \cite{ps} and 
Bo\.zejko-Picardello \cite{bp}. 
We denote by $\N_0$ the set of non-negative integers,
and by $\D$ the unit disk $\{ z\in\C : |z|<1\}$.

\begin{thm}\label{thm1}
Let $\G$ be a hyperbolic graph with bounded degree
and $d$ be the graph distance on $\G$. 
Then, there exists a constant $C$ such that the following are true. 
\begin{enumerate}
\item
For every $z\in\D$, 
the Schur multiplier $\theta_z$ on $\B(\ell_2\G)$ associated 
with the kernel 
\[
\G\times\G\ni(x,y)\mapsto z^{d(x,y)}\in\C
\] 
has (completely bounded) norm at most $C|1-z|/(1-|z|)$. 
Moreover, $z\mapsto\theta_z$ is a holomorphic map from $\D$ into
the space $V_2(\G)$ of Schur multipliers. 
\item
For every $n\in\N_0$, the Schur multiplier on $\B(\ell_2\G)$ associated 
with the characteristic function of the subset 
\[
\{ (x,y)\in\G\times\G : d(x,y)=n\}
\] 
has (completely bounded) norm at most $C(n+1)$.
\item
There exists a sequence of finitely supported functions 
$f_n\colon\N_0\to[0,1]$ such that $f_n\to1$ pointwise and that 
the Schur multiplier on $\B(\ell_2\G)$ associated with the kernel 
\[
\G\times\G\ni(x,y)\mapsto f_n(d(x,y))\in[0,1]
\]
has (completely bounded) norm at most $C$ for every $n$. 
\end{enumerate}
\end{thm}
Let $\G$ be a hyperbolic group and $d$ be the word length 
distance associated with a fixed finite generating subset of $\G$. 
Then, for the sequence $f_n$ as above, 
the sequence of functions $\p_n(x)=f_n(d(e,x))$ 
satisfy the properties required for weak amenability. 
Thus we obtain the following as a corollary.
\begin{thm}
Every hyperbolic group is weakly amenable. 
\end{thm}

This solves affirmatively a problem raised 
by Roe at the end of \cite{roe}.
We close the introduction with a few problems and remarks. 
Is it possible to construct a family of uniformly bounded 
representations as it is done in \cite{dooley,ps}? 
Is it true that a group which is hyperbolic relative to 
weakly amenable groups is again weakly amenable? 
There is no serious difficulty in extending Theorem \ref{thm1} to 
(uniformly) fine hyperbolic graphs in the sense of Bowditch \cite{bowditch}. 
Ricard and Xu \cite{rx} proved that weak amenability with constant 
one is closed under free products with finite amalgamation. 
The author is grateful to Professor Masaki Izumi for 
conversations and encouragement.
\section{Preliminary on Schur multipliers}\label{sec:schur}
Let $\G$ be a set and denote by $\B(\ell_2\G)$ the 
Banach space of bounded linear operators on $\ell_2\G$. 
We view an element $A\in\B(\ell_2\G)$ as a $\G\times\G$-matrix: 
$A=[A_{x,y}]_{x,y\in\G}$ with $A_{x,y}=\ip{A\delta_y,\delta_x}$. 
For a kernel $k\colon\G\times\G\to\C$, the Schur multiplier 
associated with $k$ is the map $m_k$ on $\B(\ell_2\G)$ 
defined by $m_k(A)=[k(x,y)A_{x,y}]$. 
We recall the necessary and sufficient condition for $m_k$ 
to be bounded (and everywhere-defined). 
See \cite{bo,pisier} for more information of completely bounded maps 
and the proof of the following theorem. 

\begin{thm}\label{schur}
Let a kernel $k\colon\G\times\G\to\C$ and a constant $C\geq0$ be given. 
Then the following are equivalent. 
\begin{enumerate}
\item 
The Schur multiplier $m_k$ is bounded and $\|m_k\|\le C$. 
\item 
The Schur multiplier $m_k$ is completely bounded and $\|m_k\|_{\cb}\le C$. 
\item 
There exist a Hilbert space $\hh$ and vectors $\zeta^+(x)$, $\zeta^-(y)$ 
in $\hh$ with norms at most $\sqrt{C}$ such that 
$\ip{\zeta^-(y),\zeta^+(x)}=k(x,y)$ for every $x,y\in\G$.
\end{enumerate}
\end{thm}
We denote by $V_2(\G)=\{ m_k : \|m_k\|<\infty\}$ the Banach space of Schur multipliers.
The above theorem says that the sesquilinear form
\[
\ell_\infty(\G,\hh)\times\ell_\infty(\G,\hh)\ni(\zeta^-,\zeta^+)
\mapsto m_k\in V_2(\G),
\]
where $k(x,y)=\ip{\zeta^-(y),\zeta^+(x)}$, is contractive for any Hilbert space $\hh$.

Let $\PP_f(\G)$ be the set of finite subsets of $\G$. 
We note that the empty set $\emptyset$ belongs to $\PP_f(\G)$. 
For $S\in\PP_f(\G)$, we define 
$\tilde{\xi}^+_S$ and $\tilde{\xi}^-_S\in\ell_2(\PP_f(\G))$ 
by 
\[
\tilde{\xi}^+_S(\omega)=\left\{\begin{array}{cl} 
1 &\mbox{ if $\omega\subset S$}\\ 
0&\mbox{ otherwise}\end{array}\right.
\mbox{ and }\
\tilde{\xi}^-_S(\omega)=\left\{\begin{array}{cl} 
(-1)^{|\omega|} &\mbox{ if $\omega\subset S$}\\ 
0&\mbox{ otherwise}\end{array}\right..
\]
We also set $\xi^+_S=\tilde{\xi}^+_S-\delta_{\emptyset}$
and $\xi^-_S=-(\tilde{\xi}^-_S-\delta_{\emptyset})$.
Note that $\xi^\pm_S\perp\xi^\pm_T$ if $S\cap T=\emptyset$. 
The following lemma is a trivial consequence of the binomial theorem. 

\begin{lem}\label{binomial}
One has 
$\|\xi^\pm_S\|^2+1=\|\tilde{\xi}^\pm_S\|^2=2^{|S|}$ and 
\[
\ip{\xi^-_T,\xi^+_S}
=1-\ip{\tilde{\xi}^-_T,\tilde{\xi}^+_S}
=\left\{\begin{array}{cl}
1 & \mbox{ if $S\cap T\neq\emptyset$}\\
0 & \mbox{ otherwise}
\end{array}\right.
\]
for every $S,T\in\PP_f(\G)$. 
\end{lem}
\section{Preliminary on hyperbolic graphs}
We recall and prove some facts of hyperbolic graphs. 
We identify a graph $\G$ with its vertex set and 
equip it with the graph distance: 
\[
d(x,y)=\min\{ n : \exists x=x_0,x_1,\ldots,x_n=y
\mbox{ such that $x_i$ and $x_{i+1}$ are adjacent}\}. 
\]
We assume the graph $\G$ to be connected so that $d$ is well-defined. 
For a subset $E\subset\G$ and $R>0$, we define the 
$R$-neighborhood of $E$ by 
\[
\NN_R(E)=\{ x\in\G : d(x,E)<R\}, 
\]
where $d(x,E)=\inf\{ d(x,y) : y\in E\}$. 
We write $B_R(x)=\NN_R(\{x\})$ for the ball with center $x$ and radius $R$.  
A geodesic path $\pathp$ is a finite or infinite sequence 
of points in $\G$ such that $d(\pathp(m),\pathp(n))=|m-n|$ 
for every $m,n$. 
Most of the time, we view a geodesic path $\pathp$ as a subset of $\G$.
We note the following fact (see e.g.,\ Lemma E.8 in \cite{bo}).
\begin{lem}
Let $\G$ be a connected graph.
Then, for any infinite geodesic path $\pathp\colon\N_0\to\G$ and
any $x\in\G$, there exists an infinite geodesic path $\pathp_x$
which starts at $x$ and eventually flows into $\pathp$ 
(i.e., the symmetric difference $\pathp\bigtriangleup\pathp_x$ is finite).
\end{lem}

\begin{defn}
We say a graph $\G$ is hyperbolic if there exists a constant $\delta>0$ 
such that for every geodesic triangle each edge is contained in 
the $\delta$-neighborhood of the union of the other two. 
We say a finitely generated group $\G$ is hyperbolic 
if its Cayley graph is hyperbolic.
Hyperbolicity is a property of $\G$ which is independent 
of the choice of the finite generating subset \cite{gdh,gromov}.
\end{defn}
From now on, we consider a hyperbolic graph $\G$ which 
has bounded degree: $\sup_x |B_R(x)|<\infty$ for every $R>0$. 
We fix $\delta>1$ satisfying the above definition.
We fix once for all an infinite geodesic path $\pathp\colon\N_0\to\G$ 
and, for every $x\in\G$, choose an infinite geodesic path 
$\pathp_x$ which starts at $x$ and eventually flows into $\pathp$.
For $x,y,w\in\G$, the Gromov product is defined by 
\[
\ip{x,y}_w=\frac{1}{2}(d(x,w)+d(y,w)-d(x,y))\geq0.
\]

See \cite{bo,gdh,gromov} for more information on hyperbolic spaces 
and the proof of the following lemma which says every geodesic 
triangle is ``thin''. 
\begin{lem}[Proposition 2.21 in \cite{gdh}]\label{thinness}
Let $x,y,w\in\G$ be arbitrary. 
Then, for any geodesic path $[x,y]$ connecting $x$ to $y$, 
one has $d(w,[x,y])\le\ip{x,y}_w+10\delta$. 
\end{lem}
\begin{lem}\label{cardinality}
For $x\in\G$ and $k\in\Z$, we set 
\[
T(x,k)=\{ w\in\NN_{100\delta}(\pathp_x) : d(w,x)\in\{k-1,k\}\,\}, 
\]
where $T(x,k)=\emptyset$ if $k<0$. 
Then, there exists a constant $R_0$ satisfying the following: 
For every $x\in\G$ and $k\in\N_0$, if we 
denote by $v$ the point on $\pathp_x$ such that $d(v,x)=k$, 
then 
\[
T(x,k)\subset B_{R_0}(v). 
\]
\end{lem}
\begin{proof}
Let $w\in T(x,k)$ and choose a point $w'$ on $\pathp_x$ 
such that $d(w,w')<100\delta$. 
Then, one has $|d(w',x)-d(w,x)|<100\delta$ and 
\[
d(w,v)\le d(w,w')+d(w',v)\le 100\delta+|d(w',x)-k|< 200\delta+1. 
\]
Thus the assertion holds for $R_0=200\delta+1$. 
\end{proof}
\begin{lem}\label{nonemptyintersection}
For $k,l\in\Z$, we set 
\[
W(k,l)=\{ (x,y)\in\G\times\G : T(x,k)\cap T(y,l)\neq\emptyset\}.
\]
Then, for every $n\in\N_0$, one has 
\[
E(n):=\{(x,y)\in\G\times\G : d(x,y)\le n\}=\bigcup_{k=0}^n W(k,n-k). 
\]
Moreover, there exists a constant $R_1$ such that 
\[
W(k,l)\cap W(k+j,l-j)=\emptyset
\] 
for all $j>R_1$.  
\end{lem}
\begin{proof}
First, if $(x,y)\in W(k,n-k)$, then one can find 
$w\in T(x,k)\cap T(y,n-k)$ and $d(x,y)\le d(x,w)+d(w,y)\le n$. 
This proves that the right hand side is contained in the left hand side. 
To prove the other inclusion, let $(x,y)$ and $n\geq d(x,y)$ be given.
Choose a point $p$ on $\pathp_x\cap\pathp_y$ 
such that $d(p,x)+d(p,y)\geq n$, 
and a geodesic path $[x,y]$ connecting $x$ to $y$. 
By Lemma~\ref{thinness}, there is a point $a$ on $[x,y]$ 
such that $d(a,p)\le\ip{x,y}_{p}+10\delta$. 
It follows that 
\[
\ip{x,p}_a+\ip{y,p}_a=d(a,p)-\ip{x,y}_p\le 10\delta.
\]
We choose a geodesic path $[a,p]$ connecting 
$a$ to $p$ and denote by $w(m)$ the point 
on $[a,p]$ such that $d(w(m),a)=m$. 
Consider the function $f(m)=d(w(m),x)+d(w(m),y)$. 
Then, one has that $f(0)=d(x,y)\le n\le d(p,x)+d(p,y)=f(d(a,p))$ 
and that $f(m+1)\le f(m)+2$ for every $m$. Therefore, 
there is $m_0\in\N_0$ such that $f(m_0)\in\{n-1,n\}$. 
We claim that $w:=w(m_0)\in T(x,k)\cap T(y,n-k)$ for $k=d(w,x)$. 
First, note that $d(w,y)=f(m_0)-k\in\{n-k-1,n-k\}$. 
Since 
\begin{align*}
\ip{x,p}_w &\le \frac{1}{2}(d(x,a)+d(a,w)+d(p,w)-d(x,p))\\
&= \frac{1}{2}(d(x,a)+d(p,a)-d(x,p))\\
&= \ip{x,p}_a\\
&\le 10\delta,
\end{align*}
one has that $d(w,\pathp_x)\le 20\delta$ by Lemma~\ref{thinness}. 
This proves that $w\in T(x,k)$. 
One proves likewise that $w\in T(y,n-k)$. 
Therefore, $T(x,k)\cap T(y,n-k)\neq\emptyset$ and $(x,y)\in W(k,n-k)$. 

Suppose now that $(x,y)\in W(k,l)\cap W(k+j,l-j)$ exists. 
We choose $v\in T(x,k)\cap T(y,l)$ and $w\in T(x,k+j)\cap T(y,l-j)$. 
Let $v_x$ (resp.\ $w_x$) be the point on $\pathp_x$ such that 
$d(v_x,x)=k$ (resp.\ $d(w_x,x)=k+j$). 
Then, by Lemma~\ref{cardinality}, one has 
$d(v,v_x)\le R_0$ and $d(w,w_x)\le R_0$. 
We choose $v_y$, $w_y$ on $\pathp_y$ likewise for $y$. 
It follows that $d(v_x,v_y)\le2R_0$ and $d(w_x,w_y)\le2R_0$. 
Choose a point $p$ on $\pathp_x\cap\pathp_y$. 
Then, one has 
$|d(v_x,p)-d(v_y,p)|\le 2R_0$ and 
$|d(w_x,p)-d(w_y,p)|\le 2R_0$. 
On the other hand, one has $d(v_x,p)=d(w_x,p)+j$ and 
$d(v_y,p)=d(w_y,p)-j$.  
It follows that 
\[
2j=d(v_x,p)-d(w_x,p)-d(v_y,p)+d(w_y,p)\le 4R_0.
\]
This proves the second assertion for $R_1=2R_0$. 
\end{proof}
\begin{lem}\label{counting}
We set 
\[
Z(k,l)=W(k,l)\cap\bigcap_{j=1}^{R_1} W(k+j,l-j)^c.
\]
Then, for every $n\in\N_0$, one has 
\[
\chi_{E(n)}=\sum_{k=0}^n\chi_{Z(k,n-k)}.
\]
\end{lem}
\begin{proof}
We first note that Lemma~\ref{nonemptyintersection} implies 
$Z(k,l)=W(k,l)\cap\bigcap_{j=1}^\infty W(k+j,l-j)^c$ and 
$\bigcup_{k=0}^n Z(k,n-k)\subset \bigcup_{k=0}^n W(k,n-k)=E(n)$. 
It is left to show that for every $(x,y)$ and $n\geq d(x,y)$, 
there exists one and only one $k$ such that $(x,y)\in Z(k,n-k)$. 
For this, we observe that $(x,y)\in Z(k,n-k)$ if and only if 
$k$ is the largest integer that satisfies $(x,y)\in W(k,n-k)$. 
\end{proof}
\section{Proof of Theorem}
\begin{prop}\label{prop}
Let $\G$ be a hyperbolic graph with bounded degree and define 
$E(n)=\{(x,y): d(x,y)\le n\}$.
Then, there exist a constant $C_0>0$, subsets $Z(k,l)\subset\G$, 
a Hilbert space $\hh$ and vectors $\eta^+_k(x)$ and $\eta^-_l(y)$ in $\hh$ 
which satisfy the following properties: 
\begin{enumerate}
\item
$\eta^\pm_m(w)\perp\eta^\pm_{m'}(w)$ for every $w\in\G$ and $m,m'\in\N_0$ with $|m-m'|\geq2$. 
\item
$\|\eta^\pm_m(w)\|\le\sqrt{C_0}$ for every $w\in\G$ and $m\in\N_0$.
\item
$\ip{\eta^-_l(y),\eta^+_k(x)}=\chi_{Z(k,l)}(x,y)$ 
for every $x,y\in\G$ and $k,l\in\N_0$. 
\item 
$\chi_{E(n)}=\sum_{k=0}^n\chi_{Z(k,n-k)}$ for every $n\in\N_0$.
\end{enumerate}
\end{prop}
\begin{proof}
We use the same notations as in the previous sections.

Let $\hh=\ell_2(\PP_f(\G))^{\otimes(1+R_1)}$ and define 
$\eta^+_k(x)$ and $\eta^-_l(y)$ in $\hh$ by 
\[
\eta^+_k(x)=\xi^+_{T(x,k)}\otimes
\tilde{\xi}^+_{T(x,k+1)}\otimes\cdots\otimes\tilde{\xi}^+_{T(x,k+R_1)}
\]
and
\[
\eta^-_l(y)=\xi^-_{T(y,l)}\otimes
\tilde{\xi}^-_{T(y,l-1)}\otimes\cdots\otimes\tilde{\xi}^-_{T(y,l-R_1)}.
\]
If $|m-m'|\geq2$, then $T(w,m)\cap T(w,m')=\emptyset$ and 
$\xi^\pm_{T(w,m)}\perp\xi^\pm_{T(w,m')}$. 
This implies the first assertion. 
By Lemma~\ref{cardinality} and the assumption that $\G$ has bounded degree, 
one has $C_1:=\sup_{w,m}|T(w,m)|\le\sup_v|B_{R_0}(v)|<\infty$.
Now the second assertion follows from Lemma~\ref{binomial} 
with $C_0=2^{C_1(1+R_1)}$. 
Finally, by Lemma~\ref{binomial}, one has 
\[
\ip{\eta^-_l(y),\eta^+_k(x)}
=\chi_{W(k,l)}(x,y)\prod_{j=1}^{R_1}\chi_{W(k+j,l-j)^c}(x,y)
=\chi_{Z(k,l)}(x,y).
\]
This proves the third assertion. 
The fourth is nothing but Lemma~\ref{counting}.
\end{proof}
\begin{proof}[Proof of Theorem \ref{thm1}]
Take $\eta^\pm_m\in\ell_\infty(\G,\hh)$ as in Proposition~\ref{prop} and 
set $C=2C_0$. 
For every $z\in\D$, we define $\zeta^\pm_z\in\ell_\infty(\G,\hh)$
by the absolutely convergent series
\[
\zeta^+_z(x)=\overline{\sqrt{1-z}}\sum_{k=0}^\infty \overline{z}^k\eta^+_k(x)
\]
and
\[
\zeta^-_z(y)=\sqrt{1-z}\sum_{l=0}^\infty z^l\eta^-_l(y),
\]
where $\sqrt{1-z}$ denotes the principal branch of 
the square root.
The construction of $\zeta^\pm_z$ draws upon \cite{ps}.
We note that the map 
$\D\ni z\mapsto (\zeta^\pm_z(w))_w\in\ell_\infty(\G,\hh)$
is (anti-)holomorphic.
By Proposition~\ref{prop}, one has
\begin{align*}
\ip{\zeta^-_z(y),\zeta^+_z(x)}
&=(1-z)\sum_{k,l}z^{k+l}\chi_{Z(k,l)}(x,y)\\
&=(1-z)\sum_{n=0}^\infty z^n\chi_{E(n)}(x,y)\\
&=(1-z)\sum_{n=d(x,y)}^\infty z^n\\
&=z^{d(x,y)}
\end{align*}
for all $x,y\in\G$, and
\begin{align*}
\|\zeta^\pm_z(w)\|^2 &\le 2|1-z|\sum_{j=0,1} 
 \|\sum_{m=0}^\infty (z^\pm)^{2m+j}\eta^\pm_{2m+j}(w)\|^2\\
&= 2|1-z|\sum_{j=0,1}\sum_{m=0}^\infty |z|^{4m+2j}\|\eta^\pm_{2m+j}(w)\|^2\\
 &\le 2|1-z|\frac{1}{1-|z|^2}C_0\\
 &< C\frac{|1-z|}{1-|z|}
\end{align*}
for all $w\in\G$.
Therefore the Schur multiplier $\theta_z$ associated with the kernel $z^d$ has 
(completely bounded) norm at most $C|1-z|/(1-|z|)$ 
by Theorem~\ref{schur}. Moreover, the map $\D\ni z\mapsto\theta_z\in V_2(\G)$
is holomorphic.

For the second assertion, we simply write $\|Z\|$ for the 
(completely bounded) norm of the Schur multiplier associated with 
the characteristic function $\chi_Z$ of a subset $Z\subset\G\times\G$. 
By Proposition~\ref{prop} and Theorem~\ref{schur}, 
one has 
\[
\|E(n)\|\le\sum_{k=0}^n\|Z(k,n-k)\|\le C_0(n+1).
\]
and $\|\{(x,y):d(x,y)=n\}\|=\|E(n)\setminus E(n-1)\|\le C(n+1)$. 
This proves the second assertion.
The third assertion follows from the previous two, by choosing 
$f_n(d)=\chi_{E(K_n)}(d)r_n^d$ for suitable 
$0<r_n<1$ and $K_n\in\N_0$ with $r_n\to1$ and $K_n\to\infty$. 
We refer to \cite{bp,haagerup} for the proof of this fact. 
\end{proof}

\end{document}